\documentclass[11pt,a4paper]{article}
\pdfoutput=1
\RequirePackage{amsmath}%
\usepackage{amsfonts}%
\usepackage{amssymb}%
\usepackage{url}%
\usepackage{color}

\usepackage{floatflt,graphicx,graphics}
\usepackage{subfigure}

\newtheorem{example}{Example}

\usepackage{url}
\usepackage[title]{appendix}%

\let\Re\relax
\DeclareMathOperator{\Re}{\mathrm{Re}}
\let\Im\relax
\DeclareMathOperator{\Im}{\mathrm{Im}}

\renewcommand{\i}{\mathrm{i}}
\newcommand{\bI}{{\bf I}}
\newcommand{\bM}{{\bf M}}
\newcommand{\bN}{{\bf N}}
\newcommand{\CC}{{\mathbb C}}
\newcommand{\RR}{{\mathbb R}}
\newcommand{\DD}{{\mathbb D}}

\usepackage{arcs}

\begin{document}
\title{Numerical computation of Stephenson's $g$-functions in multiply connected domains}
\author{Christopher C. Green$^{\rm a}$ \& Mohamed M S Nasser$^{\rm a}$}
	
\date{}
\maketitle
 	
\vskip-0.8cm %
\centerline{$^{\rm a}$Department of Mathematics, Statistics \& Physics, Wichita State University,} %
\centerline{Wichita, KS 67260-0033, USA}%
\centerline{\tt christopher.green@wichita.edu, mms.nasser@wichita.edu}%

\begin{abstract}

There has been much recent attention on $h$-functions, so named since they describe the distribution of harmonic measure for a given multiply connected domain with respect to some basepoint. 
In this paper, we focus on a closely related function to the $h$-function, known as the $g$-function, which originally stemmed from questions posed by Stephenson in~\cite{problems}.
Computing the values of the $g$-function for a given planar domain and some basepoint in this domain requires solving a Dirichlet boundary value problem whose domain and boundary condition change depending on the input argument of the $g$-function.
We use a well-established boundary integral equation method to solve the relevant Dirichlet boundary value problems and plot various graphs of the $g$-functions for different multiply connected circular and rectilinear slit domains. 

\end{abstract}

\begin{center}
\begin{quotation}
{\noindent {{\bf Keywords}.\;\; $g$-function, multiply connected domain, conformal mapping, boundary integral equation, generalized Neumann kernel}%
}%
\end{quotation}
\end{center}


\section{Introduction}

\begin{figure}[ht] %
	\centering{
	\hfill{\includegraphics[width=0.3\textwidth]{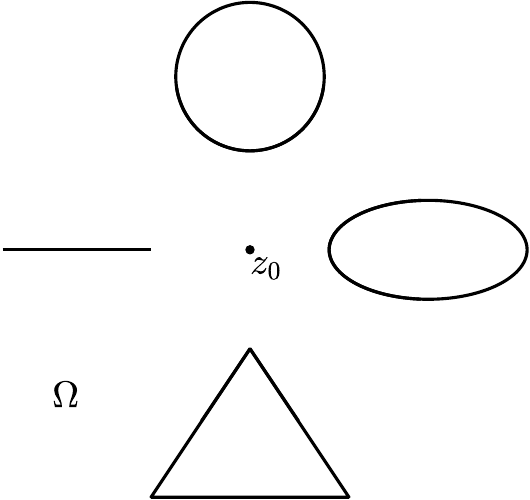}}
	\hfill
	\hfill{\includegraphics[width=0.3\textwidth]{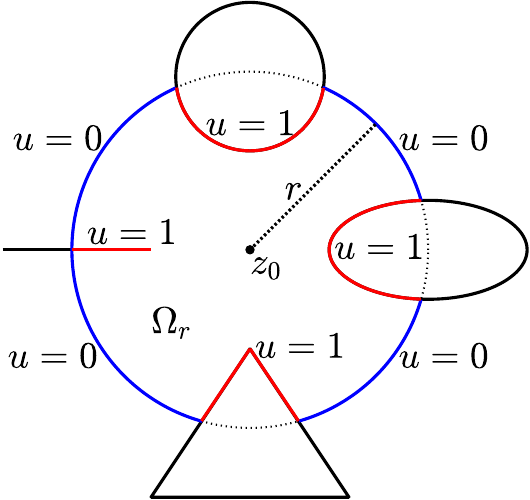}}
	\hfill
	{\includegraphics[width=0.3\textwidth]{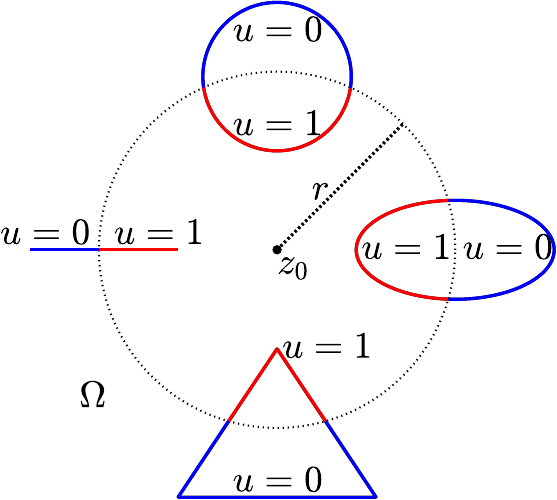}}\hfill}
\caption{Schematics of an example of an unbounded domain $\Omega$ and a basepoint $z_0\in\Omega$ (left), the domain $\Omega_r$ with the boundary conditions and the capture circle $\mathcal{C}_r$ for the $g$-function (center), and the domain $\Omega$ with the boundary conditions and the capture circle for the $h$-function (right).}
	\label{fig:gh}
\end{figure}

\begin{figure}[ht] %
	\centering{
		\hfill{\includegraphics[width=0.3\textwidth]{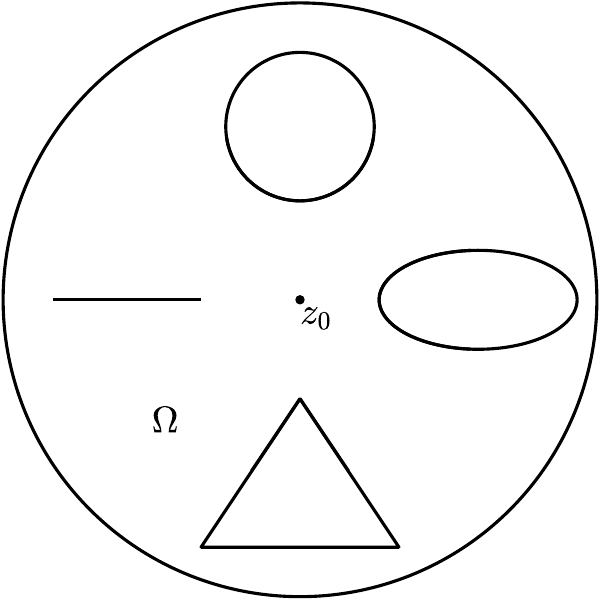}}
		\hfill
		\hfill{\includegraphics[width=0.3\textwidth]{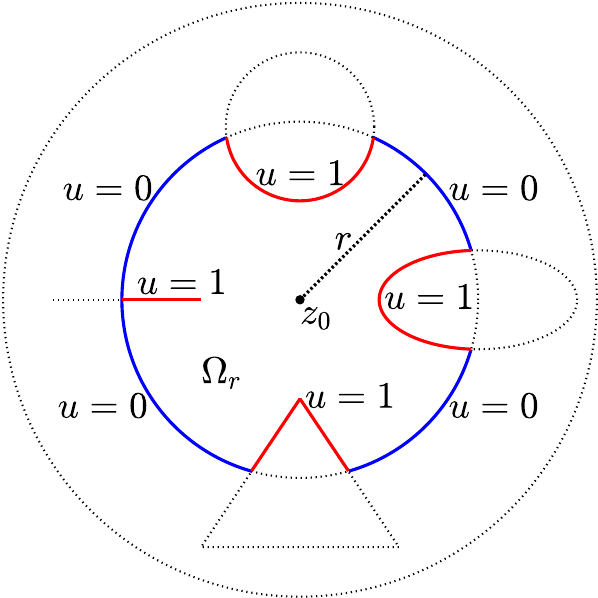}}
        \hfill
		{\includegraphics[width=0.3\textwidth]{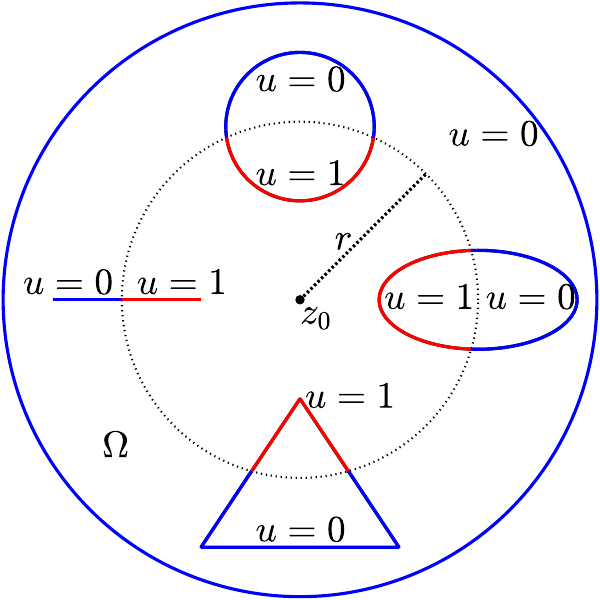}}\hfill}
	\caption{Schematics of an example of a bounded domain $\Omega$ and a basepoint $z_0\in\Omega$ (left), the domain $\Omega_r$ with the boundary conditions and the capture circle $\mathcal{C}_r$ for the $g$-function (center), and the domain $\Omega$ with the boundary conditions and the capture circle for the $h$-function (right).}
	\label{fig:ghbd}
\end{figure}

Let $\Omega$ be a domain in the extended complex plane $\overline{\CC}=\CC\cup\{\infty\}$ and let $z_0$ be a given basepoint in $\Omega$.
We assume that $\Omega$ is either an unbounded domain of connectivity $\ell$ or a bounded domain of connectivity $\ell+1$ where $\ell\ge 1$ (for example, see Figures~\ref{fig:gh} and~\ref{fig:ghbd} when $\ell=4$).  
For $r>0$, let $\Omega_r$ be the connected component of $\Omega\cap B(z_0,r)$ which contains $z_0$ and let $E_r=\partial\Omega_r\cap\overline{B(z_0,r)}$ where $B(z_0,r)$ is the open disk with center $z_0$ and radius $r$. 
Unlike the given domain $\Omega$, the domain $\Omega_r$ is not fixed and changes as $r$ increases. 
Note that the domain $\Omega_r$ is always bounded with $m+1$ boundary components, where $0\le m\le \ell$ (i.e., the domain $\Omega_r$ could be simply connected or multiply connected depending on $r$). 
We refer to $\partial B(z_0,r)$ as a `capture circle' of radius $r$ and center $z_0$ and we denote it by $\mathcal{C}_r$.

For $E_r\ne\emptyset$ and $\partial\Omega_r\backslash E_r\ne\emptyset$, the harmonic measure of $E_r$ with respect to $\Omega_r$ is the $C^2(\Omega_r)$ function $u: \Omega_r \to (0,1)$ satisfying the Laplace equation
\[
\nabla^2 u = 0
\]
in $\Omega_r$, with $u(z) \to 1$ when $z \to E_r $ and $u(z) \to 0$ when $z\to\partial\Omega_r\backslash E_r$.
Harmonic measure is a key concept in potential theory and has numerous applications to geometric function theory~\cite{AVV,CroBook,CM07,Gar,Hen,Tsu}.
The harmonic measure of $E_r$ with respect to $\Omega_r$ calculated at the point $z\in\Omega_r$ will be denoted by $\omega_{\rm hm}(z,E_r,\Omega_r)$.

\begin{figure}[ht] %
\centering
\subfigure[]{\includegraphics[width=0.24\textwidth]{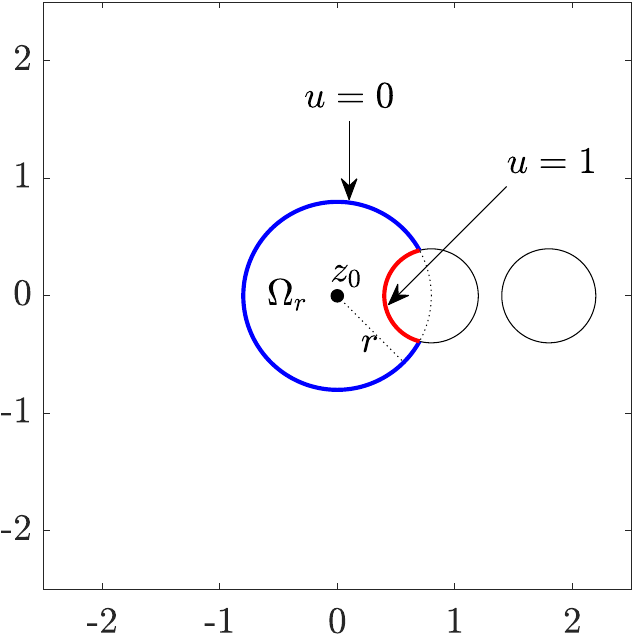}}
\subfigure[]{\includegraphics[width=0.24\textwidth]{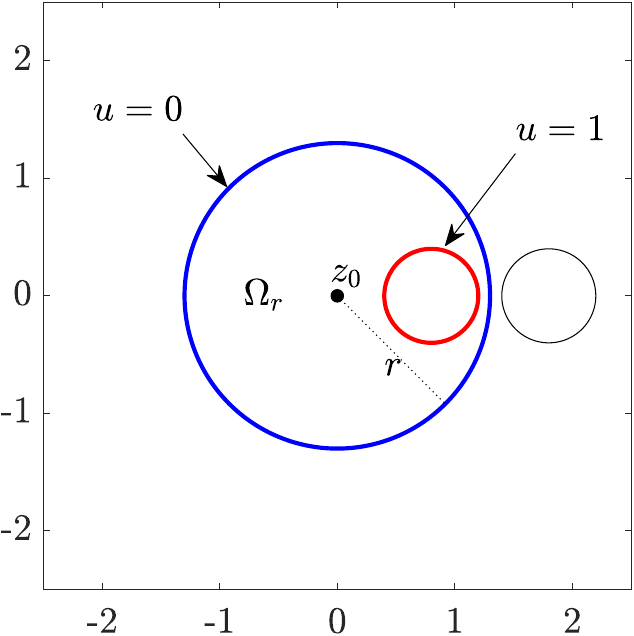}}
\subfigure[]{\includegraphics[width=0.24\textwidth]{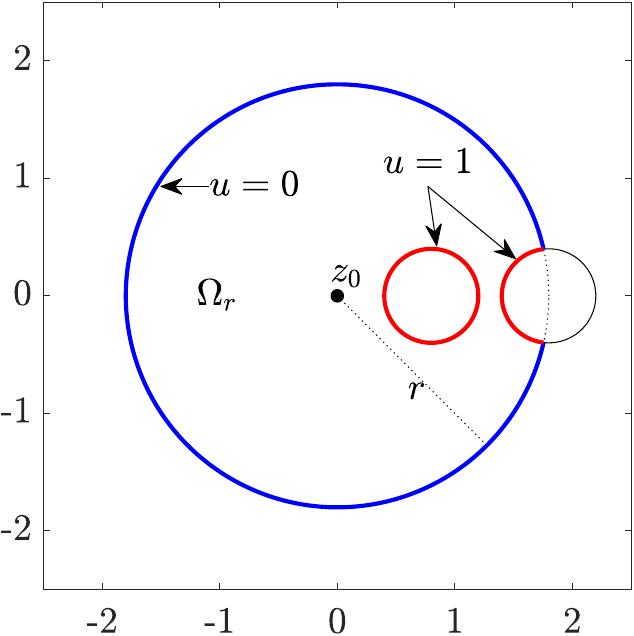}}
\subfigure[]{\includegraphics[width=0.24\textwidth]{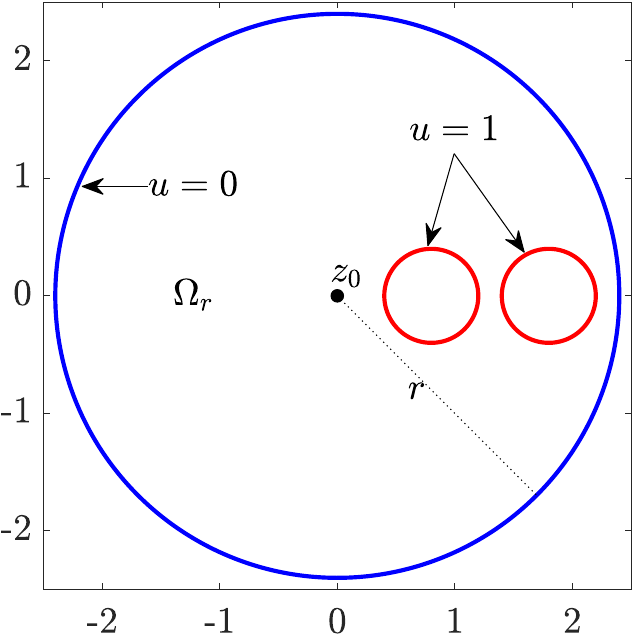}}
\caption{Schematics of the four different domain regimes $\Omega_r$ as the radius $r$ of the capture circle $\mathcal{C}_r$ increases for the doubly connected domain $\Omega$ exterior to two disks of equal radius $0.4$ centered at $0.8$ and $1.6$ with the basepoint $z_0=0$. Here, $m=0$ in (a), $m=1$ in (b) and (c), and $m=2$ in (d).}
\label{fig:g2D}
\end{figure}

\begin{figure}[ht] %
	\centering
	\subfigure[]{\includegraphics[width=0.24\textwidth]{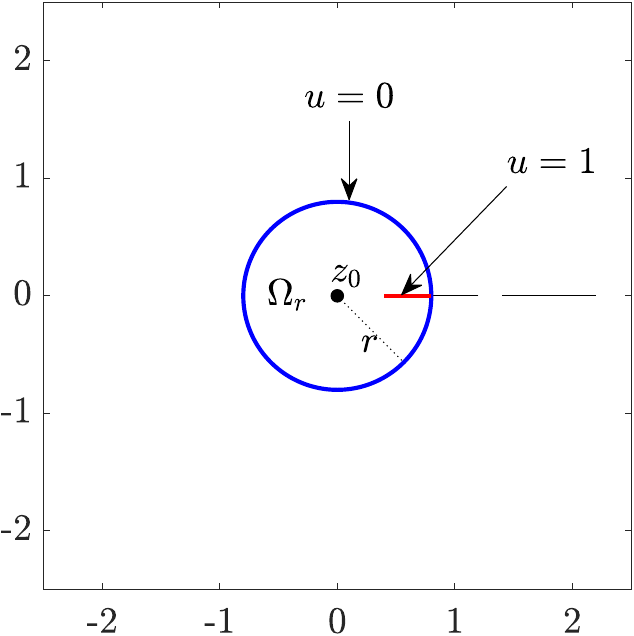}}
	\subfigure[]{\includegraphics[width=0.24\textwidth]{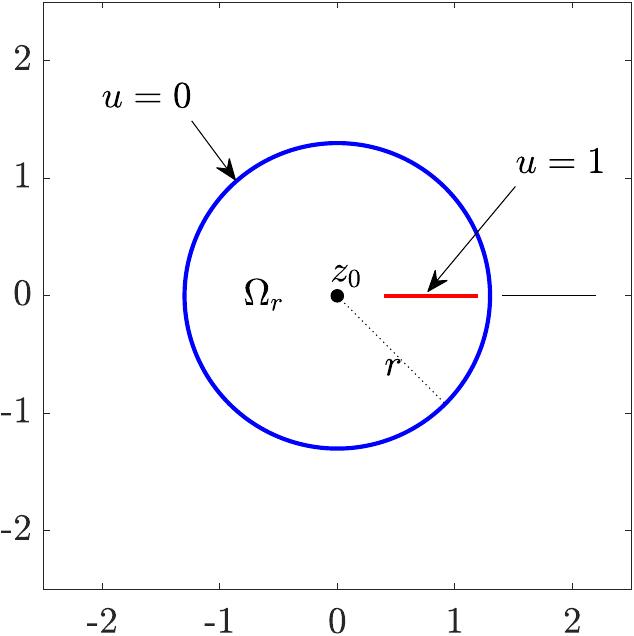}}
	\subfigure[]{\includegraphics[width=0.24\textwidth]{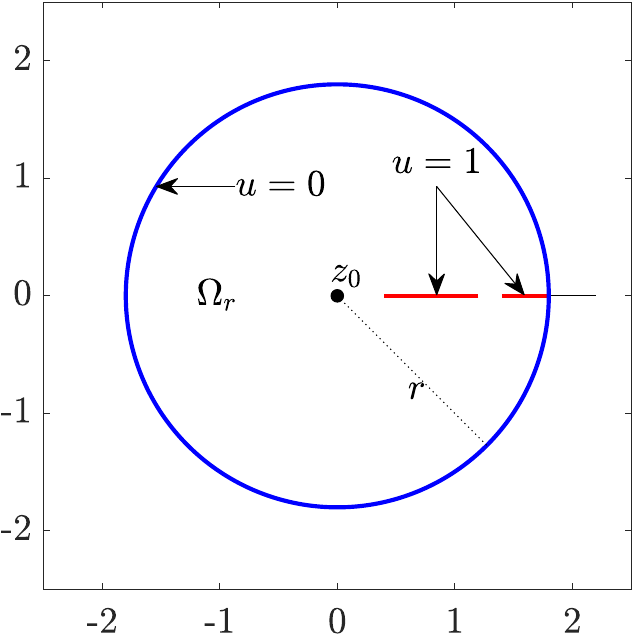}}
	\subfigure[]{\includegraphics[width=0.24\textwidth]{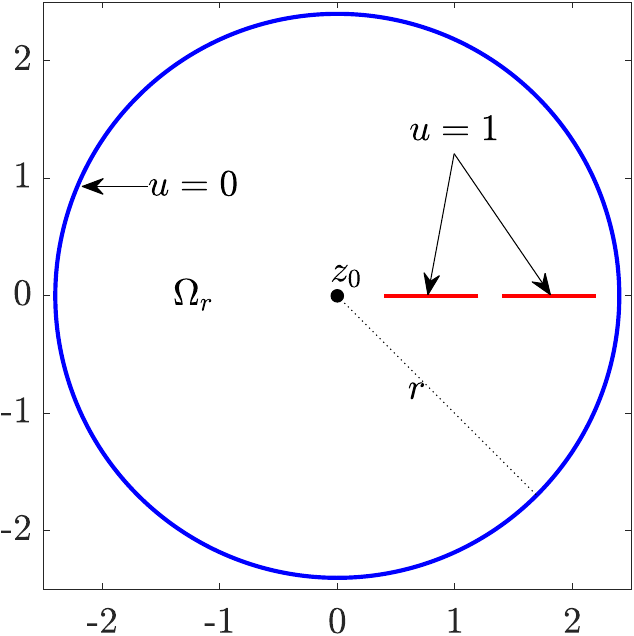}}
	\caption{Schematics of the four different domain regimes $\Omega_r$ as the radius $r$ of the capture circle $\mathcal{C}_r$ increases for the doubly connected domain $\Omega$ exterior to the two segments $[0.4,1.2]$ and $[1.4,2.2]$ of equal length with the basepoint $z_0=0$. Here, $m=0$ in (a), $m=1$ in (b) and (c), and $m=2$ in (d).}
	\label{fig:g2S}
\end{figure}

The Stephenson's $g$-function (henceforth referred to simply as the `$g$-function') associated with the domain $\Omega$ with respect to the basepoint $z_0$, $g:[0,\infty)\mapsto[0,1]$, is defined by
\begin{equation}\label{eq:g(r)-def}
g(r)=\omega_{\rm hm}(z_0,E_r,\Omega_r)
\end{equation}
where $g(0)=0$. The $g$-function, which was introduced by Stephenson in~\cite[Problem~6.116]{problems}, is a non-decreasing piecewise continuous function. 
It follows from this definition of the $g$-function that $g(r)=u(z_0)$ where $u(z)$ is the unique solution of the following Dirichlet boundary value problem (BVP):
\begin{subequations}\label{eq:bdv-u}
	\begin{align}
		\label{eq:u-LapG}
		\nabla^2 u(z) &= 0, \quad \mbox{~~~~ }z\in \Omega_r; \\
		\label{eq:u-1G}
		u(z)&= \gamma(z), \quad \mbox{ }z\in \partial \Omega_r;  
	\end{align}
\end{subequations} 
where
\begin{equation}\label{eq:bdv-gam}
	\gamma(z)= \left\{
	\begin{array}{cc} 
		1, & z\in E_r, \\ 
		0, & z\in \partial\Omega_r\backslash E_r. \\ 
	\end{array}
	\right.
\end{equation}
This is illustrated in Figure~\ref{fig:gh} for unbounded multiply connected domains $\Omega$ and in Figure~\ref{fig:ghbd} for bounded domains. See also Figures~\ref{fig:g2D} and~\ref{fig:g2S} when the domain $\Omega$ is the doubly connected domain exterior to two circles (Figure~\ref{fig:g2D}) and exterior to two slits (Figure~\ref{fig:g2S}) with the basepoint fixed to be at $z_0=0$. 

The function which is strongly associated with the $g$-function is the so-called harmonic-measure distribution function, dubbed the $h$-function, which is defined also with respect to the planar domain $\Omega$ and the basepoint $z_0$. It is equal to the value of the harmonic measure of the portion $E_r$ of the boundary with respect to $\Omega$ at $z_0$:
\[
h(r)=\omega_{\rm hm}(z_0, E_r,\Omega).
\]
The $h$-function can be computed by solving a Dirichlet boundary value problem similar to~\eqref{eq:bdv-u}. This is illustrated also in Figure~\ref{fig:gh} for unbounded multiply connected domains $\Omega$ and in Figure~\ref{fig:ghbd} for bounded domains. 
The $h$-function was also introduced by Stephenson in~\cite{problems} and was
first studied in depth by Walden \& Ward~\cite{WW}.

Naturally, owing to the connection with harmonic measure, a physical interpretation of the values of the $g$-functions -- as for the $h$-functions -- can be given in terms of a Brownian particle released into $\Omega$ from the point $z_0$ (see~\cite{Kak44,Kak45,Kar,survey} for the relation between Brownian motion and harmonic functions). More precisely, for each assignment of $r$, the value $g(r)$ is the probability that the Brownian particle will first exit the domain $\Omega_r$ through the portion $E_r$ of the boundary $\partial \Omega_r$. These hitting probabilities described by the $g$-function have particle trajectories which are confined strictly to the interior of the capture circle. On the other hand, given an $r>0$, the value $h(r)$ is the probability that a Brownian walker will first exit the domain $\Omega$ through the portion $E_r$ of the boundary $\partial \Omega$. In terms of the set of admissible trajectories of the Brownian particle, there may be trajectories which wander into the exterior of the capture circle. 
By the monotonicity of harmonic measure~\cite[p.~252]{Hen}, it can be shown that $g(r) \leq h(r)$ for a given domain $\Omega$ and a given basepoint $z_0$. 
Further, if $\Omega$ is unbounded, then the value of the $h$-function will be $1$ when the capture circle covers all the boundary components of $\Omega$. However, the $g$-function will never attain the value $1$ when $\Omega$ is unbounded even if all boundary components have been covered by the capture circle: $g(r)\to1$ as $r\to\infty$. On the other hand, for a bounded domain $\Omega$, the value of both functions will be $1$ when the capture circle covers all the boundary components of $\Omega$.
A schematic illustrating the main differences between the $g$-function and the $h$-function is given in Figure~\ref{fig:gh} when $\Omega$ is an unbounded domain  and in Figure~\ref{fig:ghbd} when $\Omega$ is a bounded domain.

The computation of the $h$-function has been the subject of several recent works.
In~\cite{green22,green25}, analytic formulas have been derived using the Schottky--Klein prime function for computing the $h$-function associated with several multiply connected circular and slit domains. In~\cite{green24}, a boundary integral equation methods has been presented for the computation of the $h$-function of a class of highly multiply connected symmetrical slit domains. 
For more information about the $h$-function and its computation, we refer the reader to~\cite{BW,green25,green24,green22,Mah-23,phd_thesis,Mah-24,Mah-ac,survey,WW}. 

Despite this body of work on the $h$-functions, the study of $g$-functions has so far received much less attention.
Analytic formulas for the $g$-function for several simply connected domains have been presented recently in~\cite{phd_thesis,Mah-in}.
In this paper, we present a boundary integral equation method for the numerical computation of the $g$-function in multiply connected domains. The proposed method is used effectively in this paper to make calculations of the $g$-functions associated with multiply connected circular and rectilinear slit domains. To the best of our knowledge, this is the first attempt to compute the $g$-function for multiply connected domains.

The layout of this paper is as follows. In Section~\ref{sc:ie}, we describe a boundary integral equation method that will be used to compute the $g$-function in this paper. In Sections~\ref{sc:cir} and~\ref{sc:slits}, we compute the $g$-function for unbounded multiply connected circular and rectilinear slit domains, respectively. In Section~\ref{sc:other}, we discuss the computation of the $g$-functions for other kinds of circular domains. Examples of $g$-functions in bounded domains are presented in Section~\ref{sc:bd}.  Finally, we make concluding remarks in Section~\ref{sc:con}.

\section{The integral equation method}
\label{sc:ie}

Let $G$ be a bounded simply or multiply connected domain whose boundary components are piecewise smooth Jordan curves. 
In this section, we outline a boundary integral equation (BIE) method for solving the following BVP:
\begin{subequations}\label{eq:bdv-uG}
	\begin{align}
		\label{eq:u-LapG}
		\nabla^2 u(z) &= 0, \quad \mbox{~~~~ }z\in G; \\
		\label{eq:u-1G}
		u(z)&= \gamma(z), \quad \mbox{ }z\in \Gamma=\partial G;  
	\end{align}
\end{subequations} 
in the domain $G$ where $\gamma$ is assumed to be a H\"older continuous function on the boundary~$\Gamma$. 
The domain $G$ is related to the domain $\Omega_r$ introduced in the previous section.
When the boundary components of $\Omega_r$ are piecewise smooth Jordan curves, we assume that $G=\Omega_r$ (see Section~\ref{sc:cir}). 
If some of the boundary components of $\Omega_r$ are slits, then we will compute numerically a conformally equivalent domain $G$ whose boundaries are piecewise smooth Jordan curves (see Section~\ref{sc:slits}). 
That is, $G$ is always assumed to be a bounded multiply connected domain of connectivity $m+1$ with piecewise smooth boundaries where $0\le m\le \ell$ (note that $G$ is simply connected when $m=0$).
The method presented in this section will be used later to solve the BVP~\eqref{eq:bdv-u} for various $g$-functions.

Let
\[
\Gamma=\partial G=\Gamma_0\cup\Gamma_1\cup\cdots\cup\Gamma_m
\]  
where  $\Gamma_0$ is the external boundary and oriented counterclockwise. The inner curves $\Gamma_1,\ldots,\Gamma_m$ are oriented clockwise. Each curve $\Gamma_k$ is parametrized by $\eta_k(t)$ for $t\in J_k=[0,2\pi]$, $k=0,1,\ldots,m$. 
We assume that $\eta_k(t)$ is twice continuously differentiable such that $\eta'_k(t)\ne0$ for $t\in J_k$.
If $\Gamma_k$ has corner points (but not cusps), we parametrize it as explained in~\cite{Nas-7}. 
Let $J$ be the disjoint union of the $m+1$ intervals $J_j=[0,2\pi]$, $j=0,1,\ldots,m$. We define a parametrization of the whole boundary $\Gamma$ on $J$ by 
\[
\eta(t)=\left\{
\begin{array}{cc} 
	\eta_0(t), & t\in J_0, \\ 
	\eta_1(t), & t\in J_1, \\
	\vdots \\
	\eta_m(t), & t\in J_m. \\ 
\end{array}
\right.
\]
With the parametrization $\eta(t)$, we define a complex function $A$ by
\begin{equation}\label{eq:A}
	A(t) = 	\eta(t)-\alpha, \quad t\in J,
\end{equation}
where $\alpha$ is a given point in the domain $G$. 
The generalized Neumann kernel $N(s,t)$ is 
defined for $(s,t)\in J\times J$ by
\begin{equation}\label{eq:N}
	N(s,t) =
	\frac{1}{\pi}\Im\left(\frac{A(s)}{A(t)}\frac{\eta'(t)}{\eta(t)-\eta(s)}\right).
\end{equation}
We define also the following kernel 
\begin{equation}\label{eq:M}
	M(s,t) =
	\frac{1}{\pi}\Re\left(\frac{A(s)}{A(t)}\frac{\eta'(t)}{\eta(t)-\eta(s)}\right),
\end{equation}
for $(s,t)\in J\times J$.
The integral operators with the kernels $N(s,t)$ and $M(s,t)$ are denoted by $\bN$ and $\bM$, respectively. 
Further details can be found in~\cite{Nas-CVEE,Nas-ETNA,Weg-Nas}.

For a given continuous function $\gamma$, the Dirichlet problem~\eqref{eq:bdv-u} has a unique solution $u(z)$ in $G$. This unique solution can be regarded as the real part of an analytic function $F(z)$ in $G$ which is not necessarily single-valued for $m>0$. However, the function $F(z)$ can be written as 
\begin{equation}\label{eq:F(z)}
F(z)=c+(z-\alpha)f(z)-\sum_{j=1}^{m}a_j\log(z-\alpha_j)
\end{equation}
where each $\alpha_j$ is a given point in the domain interior to the boundary component $\Gamma_j$, $j=1,2,\ldots,m$, and $a_1,\ldots,a_m$ are undetermined real constants~\cite{Mik}. Since we are interested in the real part of $F(z)$, we assume that $c$ is real which is also undetermined. 
It then follows that $f(z)$ satisfies the Riemann--Hilbert problem~\cite{Nas-gc,Weg-Nas}
\begin{equation}\label{eq:rhp}
\Re[A(t)f(\eta(t))]	=-c+\gamma_0(t)+\sum_{j=1}^m a_j\gamma_j(t)
\end{equation}
where $\gamma_0(t)=\gamma(t)$ and
\begin{equation}\label{eq:gamma-j}
\gamma_j(t)=\log|\eta(t)-\alpha_j|, \quad t\in J, \quad j=1,2,\ldots,m.
\end{equation}
Note that solving the Riemann--Hilbert problem~\eqref{eq:rhp} requires finding the unknown analytic function $f(z)$ as well as the $m+1$ unknown real constants $c,a_1,\ldots,a_m$ on the right-hand side of~\eqref{eq:rhp}. 

We choose $\alpha=z_0$, the basepoint, and hence
\begin{equation}\label{eq:u(z0)}
	u(z_0)=\Re[F(z_0)]=c-\sum_{j=1}^{m}a_j\log|z_0-\alpha_j|.
\end{equation}
Thus computing $u(z_0)$ requires finding only the $m+1$ real constants $c,a_1,\ldots,a_m$.
These constants will be computed using a method based on a boundary integral equation with the generalized Neumann kernel described below.

For $j=0,1,\ldots,m$, there exists a unique real function $\mu_j(t)$ and a unique piecewise constant function $\nu_j(t)$~\cite[Theorem~2]{Nas-gc}:
\[
\nu_j(t)=\nu_{k,j},\quad t\in J_k, \quad k=0,1,\ldots,m,
\]
with real constants $\nu_{0,j},\nu_{1,j},\ldots,\nu_{m,j}$, such that \begin{equation}\label{eq:f-half}
	f_j(\eta(t))=(\gamma_j(t)+\nu_j(t)+\i\mu_j(t))/A(t)
\end{equation}
are the boundary values of an analytic function $f_j(z)$ in the bounded domain $G$. The function $\mu_j$ is the unique solution of the boundary integral equation
\begin{equation}\label{eq:ie}
	(\bI-\bN)\mu_j=-\bM\gamma_j;
\end{equation}
and the function $\nu_j$ is given by
\begin{equation}\label{eq:h}
	\nu_j=[\bM\mu_j-(\bI-\bN)\gamma_j]/2.
\end{equation}

With the analytic functions $f(z)$ in~\eqref{eq:rhp} and $f_0(z),f_1(z),\ldots,f_m(z)$ in~\eqref{eq:f-half}, we define an analytic function $F(z)$ in $G$ by
\[
\tilde f(z)=f_0(z)+\sum_{j=1}^m a_jf_j(z)-f(z).
\]
The function $\tilde f(z)$ satisfies the Riemann--Hilbert problem
\[
\Re[A(t)\tilde f(\eta(t))]= 
\nu_0(t)+\sum_{j=1}^m a_j\nu_j(t)+c
\]
where the right-hand side is a piecewise constant function. It then follows from~\cite[Lemma~2]{Nas-gc} that
\[
\nu_0(t)+\sum_{j=1}^m a_j\nu_j(t)+c=0.
\]
Hence, the $m+1$ unknown real constants $c,a_1,\ldots,a_m$ are the components of the unique solution vector of the linear system
\begin{equation}\label{eqLsys}
\left[
\begin{array}{cccc}
	\nu_{0,1}      &\cdots &\nu_{0,m}   &1    \\
	\nu_{1,1}      &\cdots &\nu_{1,m}   &1    \\
	\vdots         &\ddots &\vdots         \\
	\nu_{m,1}      &\cdots &\nu_{m,m}   &1     \\ 
\end{array}
\right]
\left[\begin{array}{c}
	a_1    \\ \vdots \\ a_{m} \\  c
\end{array}\right]
= \left[\begin{array}{c}
	-\nu_{0,0} \\ -\nu_{1,0} \\  \vdots \\ -\nu_{m,0}   
\end{array}\right].
\end{equation}
The matrix of this linear system is a particular case of the matrix in the linear system in~\cite[Theorem~4]{Nas-gc} and the proof of the uniqueness of the solution of the linear system then follows.

When the capture circle does not intersect any of the boundary components (as in Figures~\ref{fig:g2D}(b,d)), then the function $\gamma$ assumes the following simple form:
\begin{equation}\label{eq:gam-01}
\gamma(t)=\left\{
\begin{array}{cc} 
	0, & t\in J_0, \\ 
	1, & t\in J_1, \\
	\vdots \\
	1, & t\in J_m. \\ 
\end{array}
\right.
\end{equation}
Hence, the function $f_0(z)$ satisfies the Riemann--Hilbert problem
\begin{equation}\label{eq:rhp-00}
\Re[A(t)f_0(\eta(t))]= 
\gamma_0(t)+\nu_0(t)
\end{equation}
where $\gamma_0(t)=\gamma(t)$. Here, the right-hand side is a piecewise constant function which implies that $\gamma_0(t)+\nu_0(t)=0$ and hence the right-hand side of the linear system~\eqref{eqLsys} will be the vector $[0,1,\ldots,1]^T$.

It is clear from~\eqref{eq:ie} and~\eqref{eq:h} that computing the $m+1$ real constants $c,a_1,\ldots,a_m$ requires solving the same integral equation with the generalized Neumann kernel~\eqref{eq:ie} but with $m+1$ different right-hand sides and to compute the function $\nu$ in~\eqref{eq:h} $m+1$ times. When the capture circle does not intersect any of the boundary components, these numbers both reduce to $m$.

In this paper, we use the MATLAB function \verb|fbie| from~\cite{Nas-ETNA} to approximate the solution of the integral equation~\eqref{eq:ie} and the function $\nu$ in~\eqref{eq:h}. In the function \verb|fbie|, the integral equation~\eqref{eq:ie} is discretized by the Nystr\"om method and the trapezoidal rule. 
This leads to an $(m+1)n\times(m+1)n$ linear system with a dense non-symmetric coefficient matrix, where $n$ is the number of nodes in the discretization of each boundary component. This linear system is then solved by the MATLAB's built-in \verb|gmres| function together with the MATLAB function \verb|zfmm2dpart| from the Fast Multipole Method (FMM) toolbox FMMLIB2D~\cite{Gre-Gim12}. This method has been used in several publications including domains with high connectivity, domains with corners, and slits domains. We will not present further details here and instead refer the readers to~\cite{green24,NasGre18,Nas-7,Nas-gc}. However, all MATLAB codes for the calculations presented in this paper are available at: \url{https://github.com/mmsnasser/gf}.

\section{Domains bounded by circles}\label{sc:cir}

Let $\Omega$ be an unbounded multiply connected domain in the exterior of $\ell$  non-overlapping circles $C_1,C_2,\ldots,C_\ell$ with centers $z_1,\ldots,z_\ell$ and radii $r_1,\ldots,r_\ell$, and let $z_0\in\Omega$ be a given basepoint. We define
\[
\delta_k=|z_0-z_k|,\quad k=1,2,\ldots,\ell.
\]
We assume that these circles are arranged such that $\delta_1<\delta_2<\cdots<\delta_\ell$.
Further, for a given $r>0$, we assume that the capture circle $\mathcal{C}_r$ intersects with at most one of the circles $C_1,C_2,\ldots,C_\ell$ (see Figure~\ref{fig:g2D} for an example when $\ell=2$).

There are two classes of boundary data to be considered as in the following two subsections.

\subsection{Continuous boundary data}\label{sec:cont}

Suppose that the capture circle $\mathcal{C}_r$ does not intersect any of the circles $C_1,\ldots,C_\ell$. 
For $m=0,1,2,\ldots,\ell$, we assume that the circles $C_1,\ldots,C_{m}$ are inside $\mathcal{C}_r$ and the circles $C_{m+1},\ldots,C_\ell$ are outside $\mathcal{C}_r$. 
That is, all circles  $C_1,\ldots,C_m$ are outside the capture circle $\mathcal{C}_r$ for $m=0$ and all circles  $C_1,\ldots,C_\ell$ are inside the capture circle $\mathcal{C}_r$ for $m=\ell$. See Figure~\ref{fig:g2D}(b,d) for $m=1$ and $2$, respectively. 

Note that $g(r)=0$ for all values of $r$ such that all the circles  $C_1,\ldots,C_\ell$ are outside the capture circle $\mathcal{C}_r$ (i.e. $m=0$).
Otherwise, if $m>0$, then the circles $C_1,\ldots,C_{m}$ are inside $\mathcal{C}_r$ and the circles $C_{m+1},\ldots,C_\ell$ are outside $\mathcal{C}_r$. 
In this case, the domain $\Omega_r$ is a bounded multiply connected domain of connectivity $m+1$ and
\[
\partial\Omega_r = C_1\cup\cdots\cup C_{m}\cup {\mathcal{C}}_r.
\]
Here, the domain $\Omega_r$ is of the type of domain $G$ considered in Section~\ref{sc:ie} with  $\Gamma_0={\mathcal{C}}_r$ and $\Gamma_k=C_k$ for $k=1,\ldots,m$.
Note that, in the BVP~\eqref{eq:bdv-u}, we have $E_r= C_1\cup\cdots\cup C_{m}$ and hence $\partial\Omega_r\backslash E_r={\mathcal{C}}_r$. Thus, $u(z)$ is the unique solution of the Dirichlet problem~\eqref{eq:bdv-uG} with the function $\gamma(t)$ as in~\eqref{eq:gam-01}, i.e., the function $u(z)$ is the harmonic measure of $C_1\cup\cdots\cup C_{m}$ with respect to the domain $\Omega_r$. 
It then follows from~\eqref{eq:g(r)-def} and~\eqref{eq:u(z0)} that 
\[
g(r)=u(z_0)=c-\sum_{j=1}^{m}a_j\log|z_0-\alpha_j|
\]
where the $m+1$ unknown real constants $c,a_1,\ldots,a_m$ are obtained by solving the linear system~\eqref{eqLsys}.

\subsection{Discontinuous boundary data}\label{sec:disc}

When the capture circle $\mathcal{C}_r$ intersects the circle $C_{m+1}$ for any $m=0,1,\ldots,\ell-1$, then the boundary data in the BVP~\eqref{eq:bdv-u} will be discontinuous.  
Since the center of the capture circle is the basepoint which is assumed to be $z_0=0$, the two intersection points of the two circles are
\begin{equation}\label{eq:x+iy}
e^{\i\arg(z_{m+1})}\left(x\pm\i y\right),
\end{equation}
where
\begin{equation}\label{eq:x+iy2}
x=\frac{r^2+|z_{m+1}|^2-R_{m+1}^2}{2|z_{m+1}|}, \quad y=\sqrt{r^2-x^2}.
\end{equation}

There are several possible different domain regimes $\Omega_r$ as the radius $r$ of the capture circle $\mathcal{C}_r$ increases. The domain $\Omega_r$ could be simply connected or multiply connected. For all possible cases, we need the following function $U_r(z)$ which will be needed to solve the BVP~\eqref{eq:bdv-u} when the boundary data is not continuous.

\subsubsection{The function $U_r$}\label{sec:Ur}
For $m=0,1,\ldots,\ell-1$, assume that $r\in(|z_{m+1}|-R_{m+1},|z_{m+1}|+R_{m+1})$ and the capture circle $\mathcal{C}_r$ intersects the circle $C_{m+1}$ at the two points $\xi_1$ and $\xi_2$ given by~\eqref{eq:x+iy}. 
Let $\mathcal{C}_r'$ be the arc of $\mathcal{C}_r$ that lies outside $C_{m+1}$ (colored blue in Figure~\ref{fig:g2D}(a,c)) and let $C_{m+1}'$  be the arc of the circle $C_{m+1}$ that lies inside $\mathcal{C}_r$ (colored red in Figure~\ref{fig:g2D}(a,c)). 
It then follows that 
\[
\hat{\mathcal{C}}_r=\mathcal{C}_r'\cup C_{m+1}'\cup\{\xi_1,\xi_2\}
\]
is a piecewise smooth Jordan curve, see Figure~\ref{fig:g2D}(a,c). We assume that $\hat{\mathcal{C}}_r$ is oriented counterclockwise. Let $\hat{\Omega}_r$ be the simply connected domain in the interior of $\hat{\mathcal{C}}_r$, and let $U_r(z)$ be the unique solution of the following Dirichlet BVP:
\begin{subequations}\label{eq:bdv-Ur}
	\begin{align}
		\label{eq:u-Lap}
		\nabla^2 U_r(z) &= 0 \quad \mbox{if }z\in \hat{\Omega}_r; \\
		\label{eq:u-1}
		U_r(z)&= 1 \quad \mbox{if }z\in C_{m+1}';\\
		\label{eq:u-0}
		U_r(z)&= 0 \quad \mbox{if }z\in \mathcal{C}_r'. 
	\end{align}
\end{subequations} 
Note that the boundary data of the BVP~\eqref{eq:bdv-Ur} on $\hat{\mathcal{C}}_r$ is not continuous.

We can find the exact solution of the BVP~\eqref{eq:bdv-Ur}. We first use the affine map
\[
\Phi(z)=e^{-\i\arg(z_{m+1})}(z-z_0)
\]
to translate and rotate the domain $\hat\Omega_r$ to obtain a domain $\tilde\Omega_r$ bounded by  $\hat{L}_r=L_1\cup L_2\cup\{\zeta_1,\zeta_2\}$ where $L_1$ is the image of $\mathcal{C}_r'$, $L_2$ is the image of $C_1'$, $\zeta_1=\Phi(\xi_1)$, and $\zeta_2=\Phi(\xi_2)$. Let $\zeta_3$ be the bisection point of the arc joining $\zeta_1$ to $\zeta_2$ so that $\zeta_1,\zeta_3,\zeta_2$ are ordered counterclockwise. Note that $\zeta_2=\overline{\zeta_1}$ and $\zeta_3$ is on the positive real line.
The M\"obius transform
\begin{equation}\label{eq:Psi-alp}
\Psi(z)=e^{-\i\alpha}\frac{z-\zeta_2}{z-\zeta_1}, \quad \alpha=\arg\left(\frac{r-\zeta_2}{\zeta_1-r}\right)
=2\tan^{-1}\left(\frac{r-x}{y}\right)
=2\tan^{-1}\left(\sqrt{\frac{r-x}{r+x}}\right),
\end{equation}
transplants the domain $\tilde\Omega_r$ onto the wedge 
\[
\{z=re^{\i\theta}\;|\; r>0,~0<\theta<\nu\},\quad \nu=\beta-\alpha, \quad \beta=\arg\left(\frac{\zeta_3-\zeta_2}{\zeta_3-\zeta_1}\right) =2\tan^{-1}\left(\frac{y}{x-\zeta_3}\right),
\]
and hence, using trigonometric identities, we have
\[
\nu=2\tan^{-1}\left(\frac{y^2-(r-x)(x-\zeta_3)}{y(r-\zeta_3)}\right)
\] 
which can be written as
\[
\nu=2\tan^{-1}\left(\frac{(r-x)(2r-d)}{yd}\right)
=2\tan^{-1}\left(\frac{2r-d}{d}\sqrt{\frac{r-x}{r+x}}\right), \quad d=r-\zeta_3.
\]
Thus, the solution $U_r(z)$ of the BVP~\eqref{eq:bdv-Ur} is given for $z\in\hat\Omega_r$ by
\begin{equation}\label{eq:Ur}
U_r(z)=\frac{1}{\pi}\Im \log\left(\Psi(\Phi(z))\right)^{\pi/\nu}
=\frac{1}{\nu}\arg \left(e^{-\i\alpha}\frac{\Phi(z)-\zeta_2}{\Phi(z)-\zeta_1}\right)
=\frac{1}{\nu}\arg \left(e^{-\i\alpha}\frac{z-\xi_2}{z-\xi_1}\right). 
\end{equation}

\begin{figure}[ht] %
\centerline{\hfill
\scalebox{0.4}{\includegraphics[trim=0 0 0 0,clip]{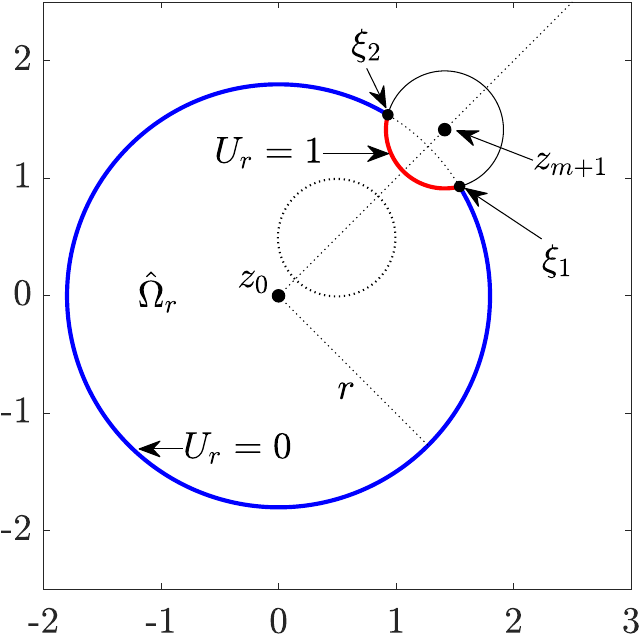}}\hfill		\scalebox{0.4}{\includegraphics[trim=0 0 0 0,clip]{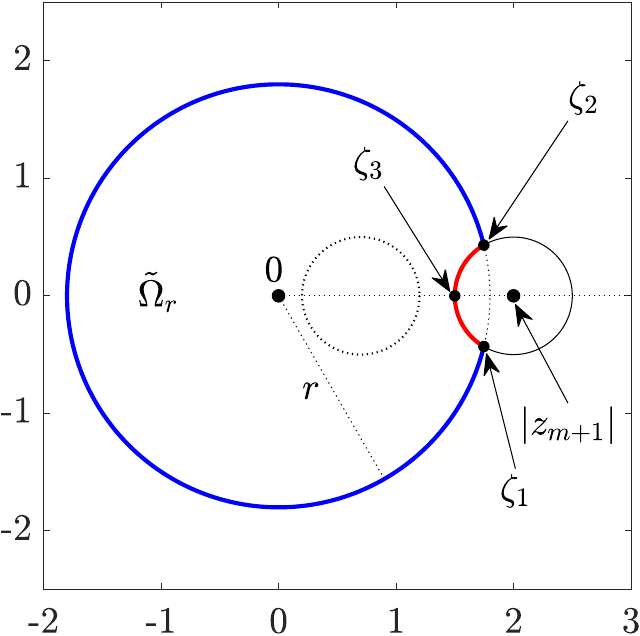}}\hfill
}
\caption{On the left, schematic of the domain $\hat{\Omega}_r$ (for $m=1$) and the boundary values of the function $U_r$ on the two arcs $C_2'$ (red) and $\mathcal{C}_r'$ (blue). On the right, the corresponding domain $\tilde\Omega_r$.}
\label{fig:g2Ur}
\end{figure}

\subsubsection{$\Omega_r$ is simply connected} 
When $m=0$ and $r\in(|z_{1}|-R_{1},|z_{1}|+R_{1})$, the capture circle $\mathcal{C}_r$ intersects the circle $C_1$ at the two points $\xi_1$ and $\xi_2$ given by~\eqref{eq:x+iy}. Note that all the other circles $C_2,\ldots,C_\ell$ are outside $\mathcal{C}_r$. 
In this case, the domain $\Omega_r$ is the same as the domain $\hat{\Omega}_r$ in Section~\ref{sec:Ur} and hence both BVPs~\eqref{eq:bdv-u} and~\eqref{eq:bdv-Ur} are identical. Thus the solution $u(z)$ to the BVP~\eqref{eq:bdv-u} is the same as the solution $U_r(z)$ to the BVP~\eqref{eq:bdv-Ur} given by~\eqref{eq:Ur}.
Since $\Phi(z_0)=0$ and $\zeta_1=\overline{\zeta_2}$, we have
\begin{equation}\label{eq:Urz01}
	g(r)=u(z_0)=U_r(z_0)=\frac{1}{\nu}\arg \left(e^{-\i\alpha}\frac{\zeta_2}{\;\overline{\zeta_2}\;}\right)
	=\frac{2\arg(\zeta_2)-\alpha}{\nu}. 
\end{equation}
It can be noted that $\arg(\zeta_2)=\tan^{-1}(y/x)$ where $x$ and $y$ are given in~\eqref{eq:x+iy2}, and hence, by~\eqref{eq:Psi-alp} and using trigonometric identities, we obtain
\[
2\arg(\zeta_2)-\alpha=2\tan^{-1}(y/x)-2\tan^{-1}\left(\frac{r-x}{y}\right)
=2\tan^{-1}\left(\frac{r-x}{y}\right)
=2\tan^{-1}\left(\sqrt{\frac{r-x}{r+x}}\right).
\]
Thus, we have
\begin{equation}\label{eq:Urz0}
	g(r)=\frac{\tan^{-1}\left(\sqrt{\frac{r-x}{r+x}}\right)}{\tan^{-1}\left(\frac{2r-d}{d}\sqrt{\frac{r-x}{r+x}}\right)}.
\end{equation}
This is equivalent to the formula given in~\cite[Eq. (21.6)]{phd_thesis}.

\subsubsection{$\Omega_r$ is multiply connected} \label{sec:cir-m}
For $1\le m\le\ell-1$ and $r\in(|z_{m+1}|-R_{m+1},|z_{m+1}|+R_{m+1})$, the capture circle $\mathcal{C}_r$ intersects the circle $C_{m+1}$ at the two points $\xi_1$ and $\xi_2$ given by~\eqref{eq:x+iy}. Note that the circles $C_1,\ldots,C_{m}$ are inside $\mathcal{C}_r$ and the circles $C_{m+2},\ldots,C_\ell$ are outside $\mathcal{C}_r$. 

In this case, the domain $\Omega_r$ is a multiply connected domain of connectivity $m+1$ and
\[
\partial\Omega_r = C_1\cup\cdots\cup C_{m}\cup \hat{\mathcal{C}}_r.
\]
The solution $u(z)$ to the BVP~\eqref{eq:bdv-u} can be written as
\begin{equation}\label{eq:u=Ur+v}
u(z)=U_r(z)+v(z)
\end{equation}
where $U_r(z)$ is given by~\eqref{eq:Ur} and $v(z)$ is a solution to the BVP
\begin{subequations}\label{eq:bdv-v}
	\begin{align}
		\label{eq:u-Lap}
		\nabla^2 v(z) &= 0 & \mbox{if }& z\in \Omega_r; \\
		\label{eq:u-1}
		v(z)&= 1-U_r(z) & \mbox{if } & z\in C_1\cup\cdots\cup C_{m}; \\
		\label{eq:u-0}
		v(z)&= 0 & \mbox{if } & z\in \hat{\mathcal{C}}_r. 
	\end{align}
\end{subequations} 
Note that the boundary conditions of the BVP~\eqref{eq:bdv-v} are now continuous and hence solving the BVP~\eqref{eq:bdv-v} is much easier than solving the original problem~\eqref{eq:bdv-u}. 
Also note that the external boundary $\hat{\mathcal{C}}_r$ of the domain $\Omega_r$ is not a circle as in the continuous data case discussed in Section~\ref{sec:cont}.

However, the domain $\Omega_r$ is still of the type of domain $G$ considered in Section~\ref{sc:ie} with  $\Gamma_0=\hat{\mathcal{C}}_r$ and $\Gamma_k=C_k$ for $k=1,\ldots,m$. It then follows from~\eqref{eq:u(z0)} that 
\[
v(z_0)=c-\sum_{j=1}^{m}a_j\log|z_0-\alpha_j|
\]
where the $m+1$ unknown real constants $c,a_1,\ldots,a_m$ are obtained by solving the linear system~\eqref{eqLsys}, and where $\gamma(\zeta)=0$ for $\zeta\in\Gamma_0$ and $\gamma(\zeta)=1-U_r(\zeta)$ for $\zeta\in\Gamma_j$ for $j=1,2,\ldots,m$. Finally, it follows from~\eqref{eq:g(r)-def} and~\eqref{eq:u=Ur+v} that
\[
g(r) = u(z_0)=U_r(z_0)+v(z_0).
\]

\subsection{Numerical examples}\label{sc:d-ex}

We assume that $\Omega$ is the unbounded multiply connected domain in the exterior of $\ell$ disks with centers $z_k$ and radii $R_k=0.4$, $k=1,2,\ldots,\ell$. We assume that the basepoint is $z_0=0$. For the centers, we consider three cases:
\begin{itemize}
	\item[(I)] $z_k=k$, $k=1,2,\ldots,\ell$.
	
	\item[(II)] $z_k=k\,e^{\theta_k\i}$, where $\theta_k=2(k-1)\pi/\ell$, $k=1,2,\ldots,\ell$.
	
	\item[(III)] $z_k=k\,e^{\tau_k\i}$, $k=1,2,\ldots,\ell$,
	where, for each $k=1,2,\ldots,\ell$, $\tau_k$ is a random number in $(0,2\pi)$.
\end{itemize}
The values of the $g$-function $g(r)$ for the three cases are shown in Figure~\ref{fig:disks} for $\ell=5$ (left) and $\ell=10$ (right).

\begin{figure}[ht] %
\centerline{\hfill
\scalebox{0.4}{\includegraphics[trim=0 0 0 0,clip]{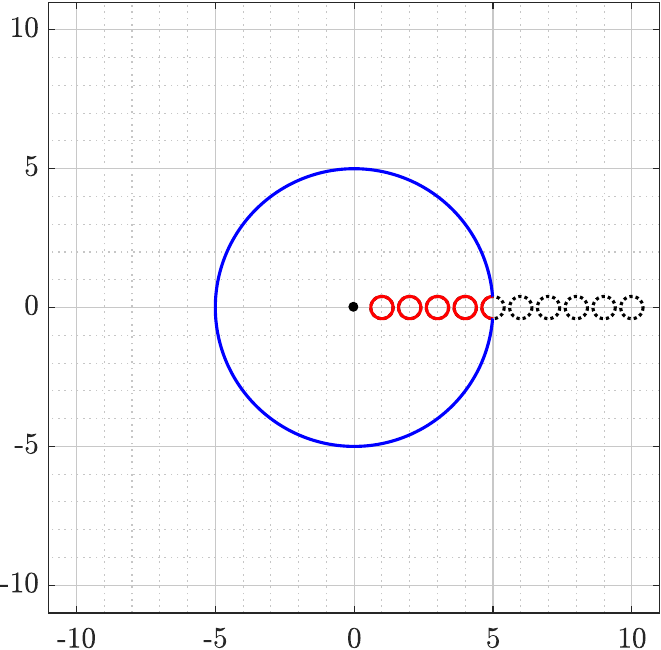}}\hfill
\scalebox{0.4}{\includegraphics[trim=0 0 0 0,clip]{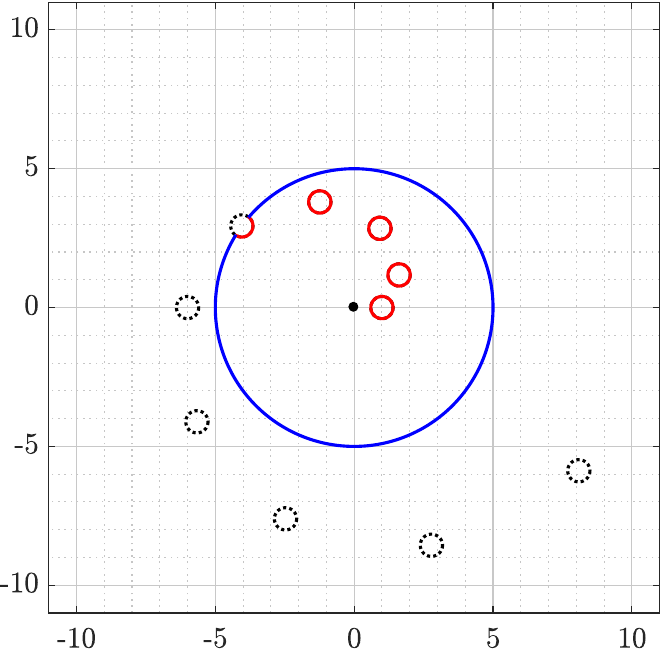}}\hfill
\scalebox{0.4}{\includegraphics[trim=0 0 0 0,clip]{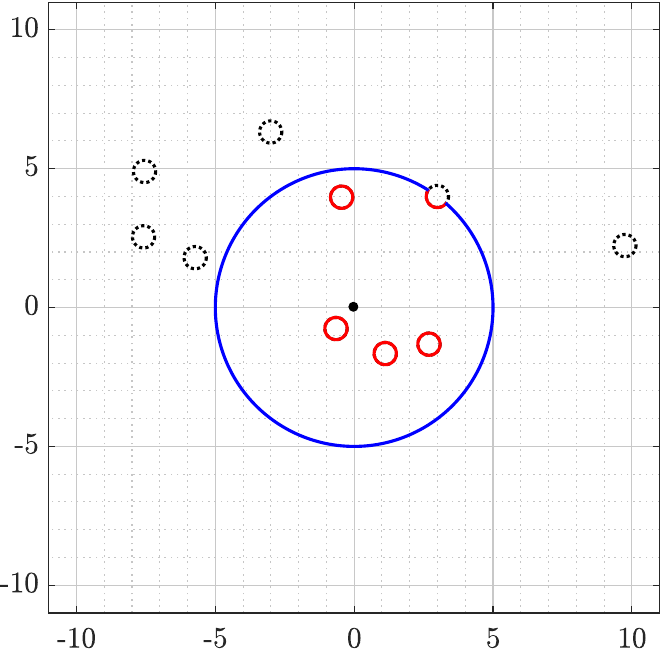}}\hfill
}
\caption{The domain $\Omega_r$ for the three unbounded multiply connected circular domain cases for $\ell=10$.}
\label{fig:mD10}
\end{figure}

\begin{figure}[ht] %
	\centerline{\hfill
		\scalebox{0.5}{\includegraphics[trim=0 0 0 0,clip]{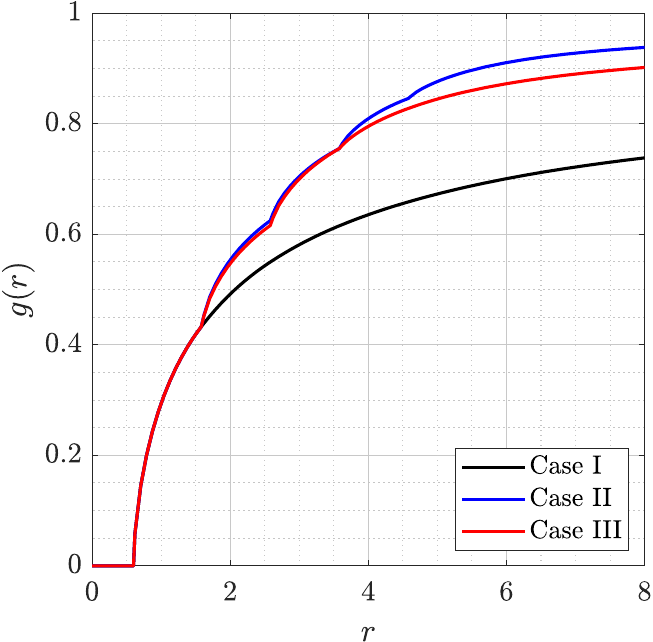}}\hfill
		\scalebox{0.5}{\includegraphics[trim=0 0 0 0,clip]{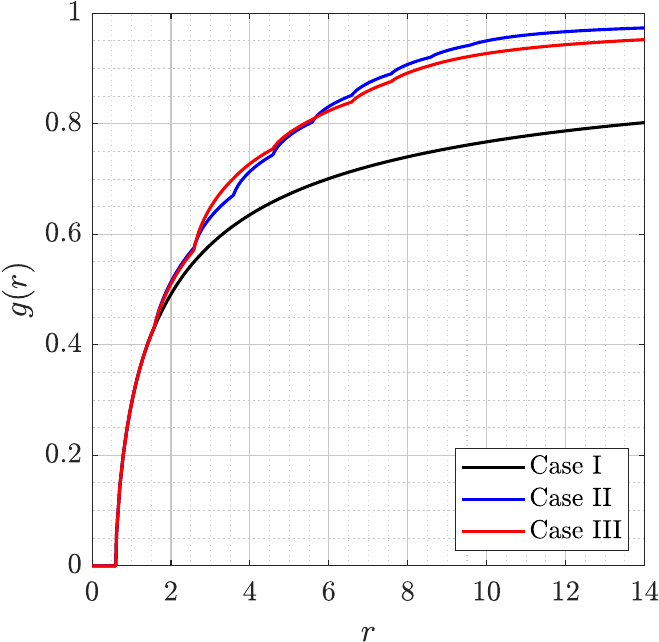}}\hfill
	}
	\caption{The values of the $g$-function $g(r)$ for $\ell=5$ (left) and $\ell=10$ (right).}
	\label{fig:disks}
\end{figure}

\section{Domains bounded by slits}\label{sc:slits}

Let $\Omega$ be an unbounded multiply connected domain in the exterior of $\ell$ non-overlapping rectilinear slits $I_j=[a_j,b_j]$ where $a_j$ and $b_j$ are real numbers. We assume that these slits are arranged such that
\[
0<a_1<b_1<a_2<b_2<\cdots<a_\ell<b_\ell.
\]
For a given $r>0$, the capture circle $\mathcal{C}_r$ intersects at most one of the slits $I_1,I_2,\ldots,I_\ell$ (see Figure~\ref{fig:g2S} for an example when $\ell=2$). 

As in the previous section, we also have here two classes of boundary data to be considered as in the following subsections.

\subsection{Continuous boundary data}\label{sec:cont2}

Suppose that the capture circle $\mathcal{C}_r$ does not intersect any of the slits $I_1,\ldots,I_\ell$. 
For $m=0,1,2,\ldots,\ell$, we assume that the slits $I_1,\ldots,I_{m}$ are inside $\mathcal{C}_r$ and the slits $I_{m+1},\ldots,I_\ell$ are outside $\mathcal{C}_r$. That is, all slits  $I_1,\ldots,I_\ell$ are outside the capture circle $\mathcal{C}_r$ for $m=0$ and all slits  $I_1,\ldots,I_\ell$ are inside the capture circle $\mathcal{C}_r$ for $m=\ell$. See Figure~\ref{fig:g2S}(b,d) for $m=1$ and $2$, respectively. 

For all values of $r$ such that all slits  $I_1,\ldots,I_\ell$ are outside the capture circle $\mathcal{C}_r$ (i.e. $m=0$), we have $g(r)=0$.
If $m\ge1$, the domain $\Omega_r$ is a bounded multiply connected domain of connectivity $m+1$ and
\[
\partial\Omega_r = I_1\cup\cdots\cup I_{m}\cup {\mathcal{C}}_r.
\]
Note that, in the BVP~\eqref{eq:bdv-u}, we have $E_r= I_1\cup\cdots\cup I_{m}$ and hence $\partial\Omega_r\backslash E_r={\mathcal{C}}_r$. 
That is, the function $u(z)$ is the harmonic measure of $I_1\cup\cdots\cup I_{m}$ with respect to the domain $\Omega_r$.

The domain $\Omega_r$ here is not bounded by Jordan curves and hence the BIE method presented in Section~\ref{sc:ie} is not directly applicable to such a domain. This obstacle will be overcome using conformal mappings. Consider the unbounded multiply connected domain $\hat \Omega$ lying in the exterior of the $m$ slits $I_1,\cdots, I_{m}$. We can use the iterative method presented in~\cite{NasGre18} to find a conformally equivalent unbounded multiply connected domain $\hat G$ lying in the exterior of $m$ circles $\Gamma_1,\cdots, \Gamma_{m}$ and a conformal mapping $w=\phi(z)$ from $\hat G$ onto $\hat\Omega$. We omit the details of the iterative method here and refer the reader to~\cite{NasGre18} (see also~\cite{green24}). Note that the circle ${\mathcal{C}}_r$, which is the outer boundary of $\Omega_r$, is within the domain $\hat\Omega$. Using the inverse mapping $z=\phi^{-1}(w)$, the circle ${\mathcal{C}}_r$ will be mapped onto a smooth Jordan curve $\Gamma_0$ surrounding the $m$ circles $\Gamma_1,\cdots, \Gamma_{m}$. Now, let $G$ be the bounded multiply connected domain in the interior of $\Gamma_0$ and in the exterior of  $\Gamma_1,\cdots, \Gamma_{m}$. Then $G$ is conformally equivalent to $\Omega_r$ and $w=\phi(z)$ is a conformal mapping from $G$ onto $\Omega_r$. Since the Dirichlet BVP is invariant under conformal mapping, then the unique solution $u(z)$ of the BVP~\eqref{eq:bdv-u} is given by $u(z)=U(\phi^{-1}(z))$ where $U(z)$ is the unique solution of the following Dirichlet BVP:
\begin{subequations}\label{eq:bdv-U}
	\begin{align}
		\label{eq:u-LapG}
		\nabla^2 U(z) &= 0, \quad \mbox{~~~~ }z\in G; \\
		\label{eq:u-1G}
		U(z)&= \hat\gamma(z), \quad \mbox{ }z\in \Gamma=\partial G;  
	\end{align}
\end{subequations} 
where
\begin{equation}\label{eq:bdv-gamU}
	\hat\gamma(z)= \left\{
	\begin{array}{cc} 
		0, & z\in \Gamma_0, \\ 
		1, & z\in \Gamma_j,\quad j=1,\ldots,m. \\ 
	\end{array}
	\right.
\end{equation}
Then, 
\[
g(r)=u(z_0)=U(\phi^{-1}(z_0)),
\]
where $U(\phi^{-1}(z_0))$ can be computed using the BIE method presented in Section~\ref{sc:ie}.

\subsection{Discontinuous boundary data}\label{sec:disc2}

When the capture circle $\mathcal{C}_r$ intersects the slit $I_{m+1}$ for any $m=0,1,\ldots,\ell-1$, then the boundary data in the BVP~\eqref{eq:bdv-u} will be discontinuous.  
Since the center of the capture circle is the basepoint which is assumed to be $z_0=0$, the single intersection point is $c_{m+1}=r$.

As before, there are several possible different domain regimes $\Omega_r$ as the radius $r$ of the capture circle $\mathcal{C}_r$ increases. The domain $\Omega_r$ can be simply or multiply connected. In all cases, we need the following function $V_r(z)$ which will be used to solve the BVP~\eqref{eq:bdv-u} when the boundary data is not continuous.

\subsubsection{The function $V_r$}\label{sec:Vr}
For $m=0,1,\ldots,\ell-1$, assume that $r\in(a_{m+1},b_{m+1})$ so that the capture circle $\mathcal{C}_r$ intersects the slit $I_{m+1}$ at the point $c_{m+1}=r$. 
Let $I_k'$ be the portion of the slit $I_k$ that lies inside $\mathcal{C}_r$ (colored red in Figure~\ref{fig:g2S}(a,c)). 
Let $\hat{\Omega}_r$ be the bounded simply connected domain such that
\[
\partial \hat{\Omega}_r=\mathcal{C}_r\cup I_k',
\]
and let $V_r(z)$ be the unique solution of the following Dirichlet BVP in $\hat{\Omega}_r$:
\begin{subequations}\label{eq:bdv-Vr}
\begin{align}
\label{eq:Vr-Lap}
\nabla^2 V_r(z) &= 0 \quad \mbox{if }z\in \hat{\Omega}_r, \\
\label{eq:Vr-1}
V_r(z)&= \gamma(z) \quad \mbox{if }z\in \partial \hat{\Omega}_r, 
\end{align}
\end{subequations} 
where
\begin{equation}\label{eq:bdv-gamV}
	\gamma(z)= \left\{
	\begin{array}{cc} 
		1, & z\in I_k', \\ 
		0, & z\in \mathcal{C}_r. \\ 
	\end{array}
	\right.
\end{equation}

The exact solution of the BVP~\eqref{eq:bdv-Vr} can be found.
To this end, consider the conformal mapping
\begin{equation}\label{eq:psi}
w=\psi(z) = \i\sqrt{\left(\frac{r-z}{r+z}\right)^2-\left(\frac{r-a_k}{r+a_k}\right)^2}
\end{equation}
which transplants $\hat{\Omega}_r$ onto the upper half-plane such that $I_k'$ is mapped onto the finite slit $[-w_k,w_k]$ on the real line and $\mathcal{C}_r$ is mapped onto $(\RR\cup\{\infty\})\backslash[-w_k,w_k]$, where
\begin{equation}\label{eq:wk}
w_k = \frac{r-a_k}{r+a_k}\in(0,1).
\end{equation}
The branch of the square root in~\eqref{eq:psi} is chosen such that $\Re\sqrt{z}\ge0$ for $z\in\CC$.
Thus, the solution $V_r(z)$ of the BVP~\eqref{eq:bdv-Vr} is given for $z\in\hat\Omega_r$ by
\begin{equation}\label{eq:Vr}
V_r(z)=\frac{1}{\pi}\arg \left(\frac{\psi(z)-w_k}{\psi(z)+w_k}\right). 
\end{equation}

\subsubsection{$\Omega_r$ is simply connected} 
When $m=0$ and $r\in(a_{1},b_{1})$, the capture circle $\mathcal{C}_r$ intersects the slit $I_{1}$ at the point $c_{1}=r$, and all the other slits $I_2,\ldots,I_\ell$ are outside $\mathcal{C}_r$. 
In this case, the domain $\Omega_r$ is the same as the domain $\hat{\Omega}_r$ in Section~\ref{sec:Vr} and hence the BVPs~\eqref{eq:bdv-u} and~\eqref{eq:bdv-Vr} are identical. Thus the solution $u(z)$ to the BVP~\eqref{eq:bdv-u} is the same as the solution $V_r(z)$ to the BVP~\eqref{eq:bdv-Vr} given by~\eqref{eq:Vr}.

It follows from~\eqref{eq:psi} that
\[
\psi(z_0)=\psi(0)=\i\sqrt{1-w_k^2}
\]
and hence~\eqref{eq:Vr} implies 
\[
V_r(z_0)
=\frac{1}{\pi}\arg \left(\frac{\i\sqrt{1-w_k^2}-w_k}{\i\sqrt{1-w_k^2}+w_k}\right)
=\frac{2}{\pi}\arg \left(\sqrt{1-w_k^2}+\i w_k\right).
\]
Finally, we have
\begin{equation}\label{eq:Vrz0}
	g(r)=V_r(z_0)
	=\frac{2}{\pi}\tan^{-1}\left(\frac{w_k}{\sqrt{1-w_k^2}}\right)
	=\frac{2}{\pi}\tan^{-1}\left(\frac{r-a_k}{2\sqrt{ra_k}}\right).
\end{equation}

\subsubsection{$\Omega_r$ is multiply connected} 
For $1\le m\le\ell-1$ and $r\in(a_{m+1},b_{m+1})$, the capture circle $\mathcal{C}_r$ intersects the slit $I_{m+1}$ at the point $c_{m+1}=r$. Note that the slits $I_1,\ldots,I_{m}$ are inside $\mathcal{C}_r$ and the slits $I_{m+2},\ldots,I_\ell$ are outside $\mathcal{C}_r$. 
In this case, the domain $\Omega_r$ is a bounded multiply connected domain of connectivity $m+1$ and
\[
\partial\Omega_r = I_1\cup\cdots\cup I_{m}\cup I_{m+1}'\cup {\mathcal{C}}_r,
\]
where $I_{m+1}'\cup {\mathcal{C}}_r$ is the outer boundary of $\Omega_r$.
We can write the solution $u(z)$ to the BVP~\eqref{eq:bdv-u} as
\[
u(z)=V_r(z)+v(z)
\]
where $V_r(z)$ is given by~\eqref{eq:Vr} and $v(z)$ is a solution to the following Dirichlet BVP:
\begin{subequations}\label{eq:bdv-v2}
\begin{align}
\label{eq:v2-Lap}
\nabla^2 v(z) &= 0 \quad \mbox{if }z\in \Omega_r;\\
\label{eq:v2-1}
v(z)&= \hat\gamma(z) \quad \mbox{if }z\in \partial\Omega_r, 
\end{align}
\end{subequations} 
where
\begin{equation}\label{eq:bdv-gamv}
\hat\gamma(z)= \left\{
	\begin{array}{cc} 
		1-V_r(z), & z\in I_1\cup\cdots\cup I_{m}, \\ 
		0, & z\in I_{m+1}'\cup {\mathcal{C}}_r. \\ 
	\end{array}
	\right.
\end{equation}

We point out that the boundary data of the BVP~\eqref{eq:bdv-v} is now continuous. However, the boundary components of the domain $\Omega_r$ are not Jordan curves and hence this problem can not be solved using the above described BIE method as in Section~\ref{sec:cir-m} when $\Omega$ was a circular domain. Thus, in the current case, we will first use conformal mappings to map $\Omega_r$ onto a domain $G$ bounded by smooth Jordan curves where we can then use our BIE method.

The mapping function $w=\psi(z)$ given by~\eqref{eq:psi} maps the external boundary $I_{m+1}'\cup {\mathcal{C}}_r$ of $\Omega_r$ onto $\tilde I_{0}=\RR\cup\{\infty\}$ and maps the slits $I_1,\cdots, I_{m}$ onto slits $\tilde I_1,\ldots,\tilde I_{m}$, respectively, on the positive imaginary axis.
Thus, $w=\psi(z)$ conformally maps the bounded multiply connected domain $\Omega_r$ onto the unbounded multiply connected domain $\tilde\Omega_r$ consists of the upper half-plane with $m$ rectilinear slits. 
For the domain $\tilde{\Omega}_r$, we use the iterative method presented in~\cite{NasGre18} to find a multiply connected circular domain $G$ in the interior of the unit circle $\Gamma_0$ and in the exterior of $m$ circles $\Gamma_1,\ldots,\Gamma_m$ and a conformal mapping $w=\Phi(\zeta)$ from the domain $G$ onto the domain  $\tilde\Omega_r$ such that $\Phi(\Gamma_j)=\tilde I_j$, $j=0,1,\ldots,m$. We can compute the circular domain $G$ such that $0\in G$, and 
\[
\Phi(\i)=\infty, \quad \Phi(0)=\psi(z_0)=\i\sqrt{1-w_k^2}
\]
where $w_k$ is given by~\eqref{eq:wk}. 
Hence, the mapping function
\[
z=\psi^{-1}(\Phi(\zeta))
\]
conformally maps the domain $G$ onto the domain $\Omega_r$ such that the inner curves $\Gamma_1,\ldots,\Gamma_{m}$ are mapped onto the slits $I_1,\ldots,I_{m}$, and the outer boundary $\Gamma_0$ is mapped onto the outer boundary $I_{m+1}'\cup {\mathcal{C}}_r$ of $\Omega_r$. Consequently, the solution $v(z)$ to the BVP~\eqref{eq:bdv-v2} is given by
\begin{equation}\label{eq:v-V}
v(z)=V(\Phi^{-1}(\psi(z))), \quad z\in\Omega_r,
\end{equation} 
where $V(\zeta)$ is the unique solution to the following Dirichlet BVP in the domain $G$:
\begin{subequations}\label{eq:bdv-v3}
\begin{align}
\label{eq:v3-Lap}
\nabla^2 V(\zeta) &= 0 \quad \mbox{if }\zeta\in G; \\
\label{eq:v3-1}
V(\zeta)&=\tilde\gamma(\zeta) \quad \mbox{if }\zeta\in \Gamma= \partial G, 
\end{align}
\end{subequations}  
where
\begin{equation}\label{eq:bdv-gamV2}
	\hat\gamma(\zeta)= \left\{
	\begin{array}{cc} 
		1-V_r(\psi^{-1}(\Phi(\zeta))), & \zeta\in \Gamma_1\cup\cdots\cup \Gamma_{m}, \\ 
		0, & \zeta\in \Gamma_{0}. \\ 
	\end{array}
	\right.
\end{equation}

The BVP~\eqref{eq:bdv-v3} can then be solved easily using the BIE method as described in Section~\ref{sc:ie}. By computing the solution $V(\zeta)$ to the BVP~\eqref{eq:bdv-v3}, we obtain the solution $v(z)$ to the BVP~\eqref{eq:bdv-v2} by~\eqref{eq:v-V}. Hence,
\[
g(r)=u(z_0)=V_r(z_0)+v(z_0)=V_r(z_0)+V(\Phi^{-1}(\psi(z_0))))
=V_r(z_0)+V(0).
\]

\subsection{Numerical examples}

\begin{figure}[ht] %
	\centerline{
		\hfill\scalebox{0.5}{\includegraphics[trim=0 0 0 0,clip]{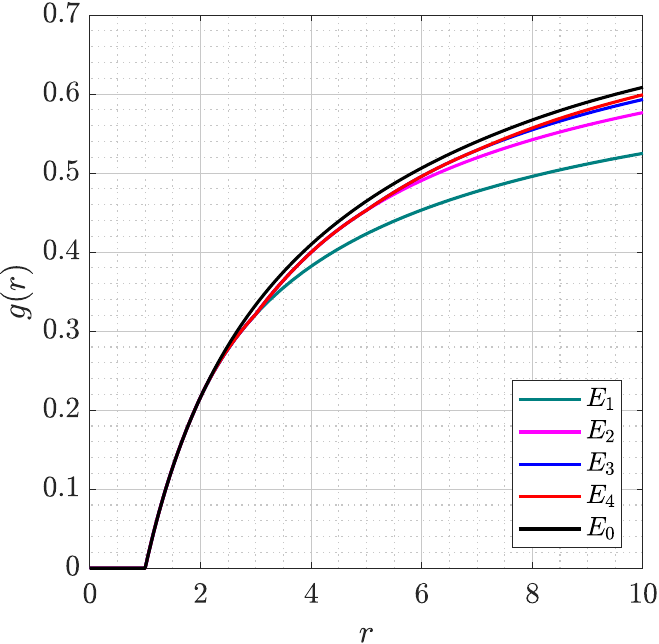}}
		\hfill
		\scalebox{0.5}{\includegraphics[trim=0 0 0 0,clip]{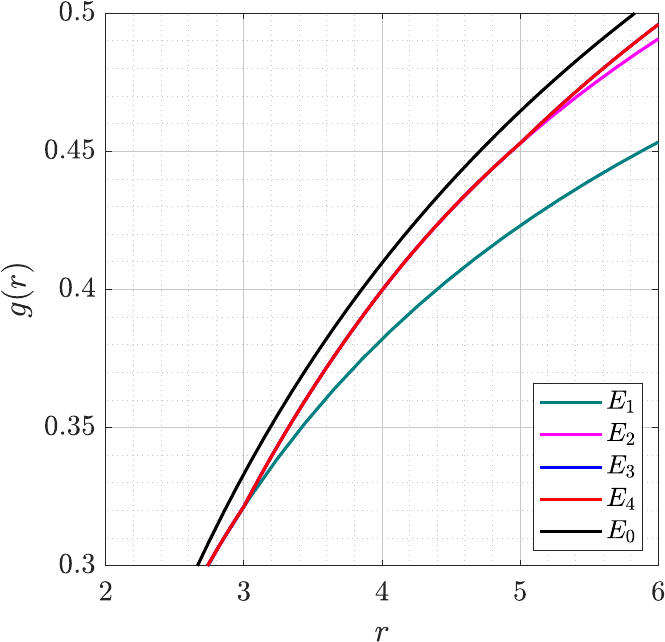}}
		\hfill
	}
	\caption{The values of the $g$-function $g(r)$ for $\Omega_k$ in Example~\ref{ex:m-int} for $k=0,1,2,3,4$ (left) and the magnification of the part of the figure on $[2,6]\times[0.3,0.5]$ (right).}
	\label{fig:IntS}
\end{figure}

\begin{example}\label{ex:m-int}{\rm
We assume that the basepoint is $z_0=0$ and $\Omega_k$ is the unbounded multiply connected domain in the exterior of the closed set $E_k$ for $k=0,1,2,3,4$ where $E_0=[1,8]$, $E_1=[1,2]$, $E_2=[1,2]\cup[3,4]$, $E_3=[1,2]\cup[3,4]\cup[5,6]$, and $E_4=[1,2]\cup[3,4]\cup[5,6]\cup[7,8]$. 
The values of the $g$-function $g(r)$, for $k=0,1,2,3,4$, computed with $n=2^{13}$, are shown in Figure~\ref{fig:IntS} (left). A magnification of the part of the figure on $[2,6]\times[0.3,0.5]$ is shown in Figure~\ref{fig:IntS} (right).

The sets $E_0$ and $E_1$ consist of only one interval and the graphs of the $g$-function $g(r)$ for these two sets only have one corner point at which $g'(r)$ is discontinuous. This occurs at $r=1$ when the capture circle first hits each of the sets $E_0$ and $E_1$. The set $E_2$ consists of two intervals and the graph of the $g$-function $g(r)$ for this set shows that $g'(r)$ is discontinuous at both $r=1$ and $r=3$ (i.e. at the left-ends of the two intervals making up the set $E_2$). Similarly, the graphs of the $g$-function $g(r)$ for the sets $E_3$ and $E_4$ show that $g'(r)$ is discontinuous at $r=1,3,5$ for $E_3$ and at $r=1,3,5,7$ for $E_4$.
}\end{example}

\begin{figure}[ht] %
\centerline{
\scalebox{0.5}{\includegraphics[trim=0 0 0 0,clip]{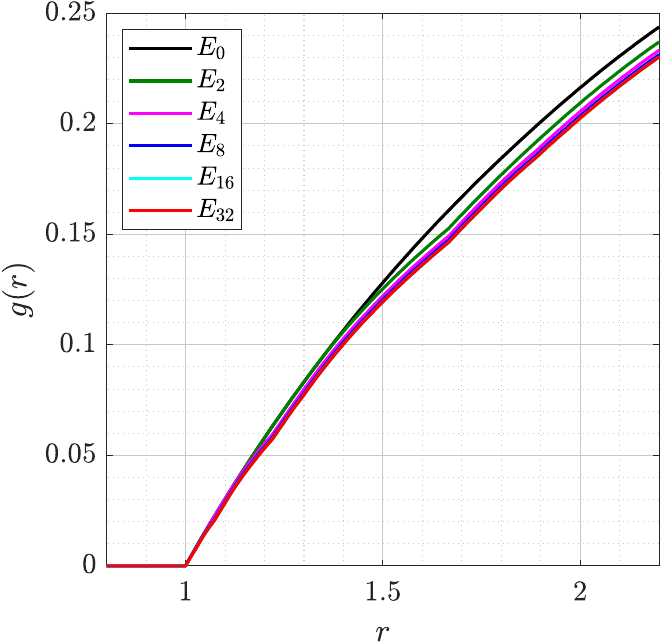}}
}
\caption{The values of the $g$-function $g(r)$ for $\Omega_k$ in Example~\ref{ex:m-cant} for $k=0,2,4,8,16,32$.}
\label{fig:CantS}
\end{figure}

\begin{example}\label{ex:m-cant}{\rm
Let the sets $E_k$ be defined recursively by
\[
E_k=\frac{1}{3}\left(E_{k-1}+2\right)\bigcup\frac{1}{3}\left(E_{k-1}+4\right), \quad k\ge1,
\]
with $E_0=[1,2]$. Note that
\[
\bigcap_{j=0}^{\infty}E_j
\]
is the middle-thirds Cantor set of the interval $[1,2]$. We assume that the basepoint is $z_0=0$ and $\Omega_k$ is the unbounded multiply connected domain in the exterior of the closed set $E_k$. The graphs of the $g$-function $g(r)$, for $k=0,2,4,8,16,32$, computed with $n=2^{13}$ are shown in Figure~\ref{fig:CantS}.	These graphs illustrate that the $g$-function for the domain $\Omega_k$ is non-constant and exhibits a finite number of points where its first derivative is discontinuous.  
These points correspond to the values of $r$ at the $2^k$ left-ends of the intervals making up the set $E_k$.
Note that, for this domain $\Omega_k$, the graphs of the $h$-function for $k=16$ and $k=32$ are presented in~\cite[Fig~8]{green24}. In contrast, these graphs consist of a finite number of horizontal line segments which correspond to those values of $r$ when the capture circle is not intersecting the slits and is otherwise an increasing function.

As is apparent in Figure~\ref{fig:CantS}, the values of the $g$-function, as a function of $k$, decrease slightly as $k$ increases. However, it seems that these graphs approach some curve as $k$ increases and it is difficult to tell the graphs apart from each other, e.g., it is not possible in Figure~\ref{fig:CantS} to distinguish between the graphs for $k=16$ and $k=32$.
}\end{example}

\section{Other domains}\label{sc:other}

\begin{figure}[ht] %
	\centerline{\hfill
		\scalebox{0.5}{\includegraphics[trim=0 0 0 0,clip]{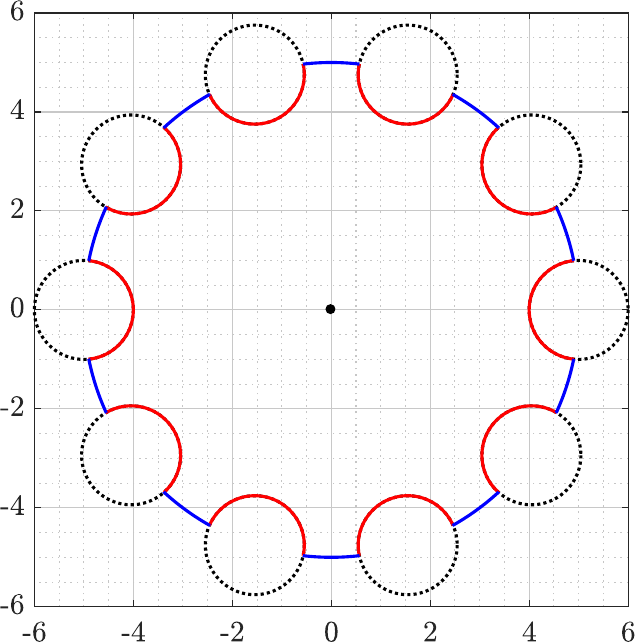}}\hfill
		\scalebox{0.5}{\includegraphics[trim=0 0 0 0,clip]{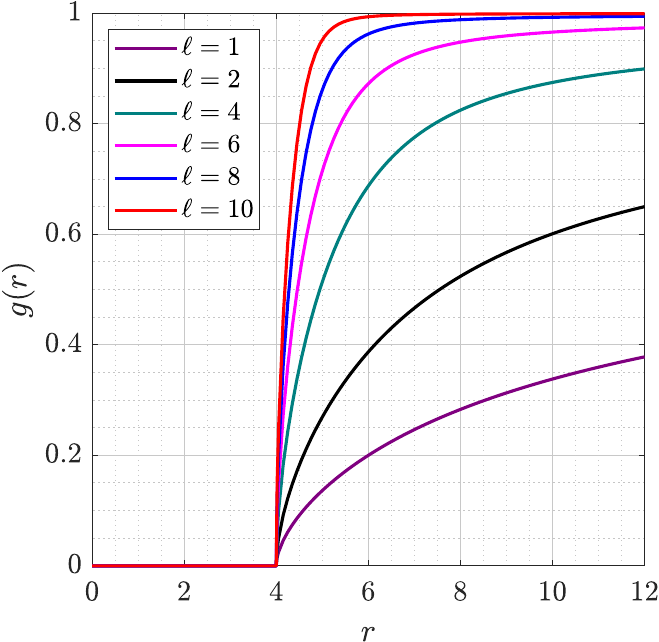}}\hfill
	}
	\caption{The domain $\Omega_r$ for the centers and radii given in~\eqref{eq:other} for $\ell=10$ and $r=5$ (left) and the values of the $g$-function $g(r)$ for several values of $\ell$ (right).}
	\label{fig:disksR}
\end{figure}

In the previous two sections, we assumed that the capture circle intersects at most one of the boundary components. The presented method can be extended to the cases when the capture circle intersects more than one boundary component. By way of example, consider the unbounded multiply connected domain $\Omega$ in the exterior of $m$ disks having centers and radii
\begin{equation}\label{eq:other}
	z_k=5\,e^{2\pi (k-1)\i/\ell}, \quad R_k=1, \quad k=1,2,\ldots,\ell.
\end{equation}
We assume that the basepoint is $z_0=0$. It is clear that $g(r)=0$ for $r\in[0,4]$. For $r\in(4,6)$, the domain $\Omega_r$ is a simply connected polycircular arc domain as shown in Figure~\ref{fig:disksR} (left) for $\ell=10$. Then
\[
g(r)=u(z_0)
\] 
where $u(z)$ is the unique solution to the BVP~\eqref{eq:bdv-u} where $E_r$ is the union of the parts of the circles $C_1,\ldots,C_\ell$ that lie inside the capture circle $\mathcal{C}_r$ (shown in red in Figure~\ref{fig:disksR} (left)). In this case, the capture circle intersects all the circles $C_1,\ldots,C_\ell$. The intersection points are denoted by $\xi_{2j-1}$ and $\xi_{2j}$, $j=1,2,\ldots,\ell$. We assume $\partial\Omega_r$ is oriented counterclockwise and that the points $\xi_1,\ldots,\xi_{2\ell}$ are ordered counterclockwise on the capture circle $\mathcal{C}_r$. 

When $r\in(4,6)$, we let $\zeta=\psi(z)$ be the conformal mapping from the simply connected polycircular arc domain $\Omega_r$ onto the unit disk $\DD$ such that $\psi(0)=0$. Then $\zeta=\psi(z)$ maps the boundary $\partial\Omega_r$ onto the unit circle $L$ such that the points $\zeta_j=\psi(\xi_j)$, $j=1,\ldots,2\ell$, are on the unit circle and oriented counterclockwise. For $k=1,2,\ldots,\ell$, let $L_k$ be the arc on the unit circle that connects $\zeta_{2k-1}$ to $\zeta_{2k}$ (in the counterclockwise direction). Let also $\hat L_k$ be the arc on the unit circle that connects $\zeta_{2k}$ to $\zeta_{2k+1}$ (in the counterclockwise direction), where $\zeta_{2\ell+1}=\zeta_1$. The mapping function $\zeta=\psi(z)$ can be computed as described in~\cite{Nas-7}.

Let $w=\phi(\zeta)$ be the M\"obius transformation from the unit disk $\DD$ onto the upper half-plane such that $\phi(\zeta_1)=\infty$ and $\phi(\zeta_{\ell+1})=0$. Hence, 
\[
w=\phi(\zeta)=\frac{\zeta^\ast-\zeta_1}{\zeta^\ast-\zeta_{\ell+1}}\frac{\zeta-\zeta_{\ell+1}}{\zeta-\zeta_1}
\]
where $\zeta^\ast=\zeta_{\lfloor(3\ell/2)+1\rfloor}$ and $\lfloor\cdot\rfloor$ denotes the floor function. The point $\zeta^*$ is chosen such that it will be mapped onto $1$ on the positive real line so that the arc connecting $\zeta_1$ and $\zeta_{\ell+1}$ (in the counterclockwise direction) is mapped onto the negative real line. Then, for $j=2,\ldots,2\ell$, the points $x_j=\phi(\zeta_j)$ will be on the real line such that $x_j<x_{j+1}$. Let $I_j$ and $\hat I_j$ be the images of the arcs $L_j$ and $\hat L_j$ under the M\"obius transformation $w=\phi(\zeta)$, respectively. Then $I_2,\ldots,I_\ell,\hat I_1,\ldots\hat I_{\ell-1}$ are finite intervals and $I_1,\hat I_\ell$ are infinite intervals on the real line.

Now, the function
\[
U(w)= \sum_{j=2}^{2\ell} \frac{(-1)^j}{\pi}\arg(w-x_j)
\]
is harmonic on the upper half-plane with 
\[
U(w)= \left\{
\begin{array}{cc} 
	1, & \zeta\in I_1\cup\cdots\cup I_{\ell}, \\ 
	0, & \zeta\in \hat I_1\cup\cdots\cup \hat I_{\ell}. \\ 
\end{array}
\right.
\]
The branch of $\arg(\cdot)$ is chosen such that $\arg(-1)=\pi$.
Then the unique solution $u(z)$ to the BVP~\eqref{eq:bdv-u} is given by
\[
u(z)=U(\phi(\psi(z)))
= \sum_{j=2}^{2\ell} \frac{(-1)^j}{\pi}\arg\left(\phi(\psi(z))-\phi(\psi(\xi_j))\right),
\]
and hence the $g$-function $g(r)$ can be computed via
\[
g(r)=u(z_0)=u(0).
\]

Finally, when $r>6$, the domain $\Omega_r$ is a bounded multiply connected circular domain of connectivity $\ell+1$. The domain $\Omega_r$ in this case is of the type considered in Section~\ref{sc:ie} with $G=\Omega_r$, $m=\ell$, and the function $\gamma$ as in~\eqref{eq:gam-01}. Hence, $g(r)=u(z_0)$ where $u(z_0)$ can be computed as in~\eqref{eq:u(z0)}.

The graphs of the $g$-function $g(r)$ for several values of $\ell$ are shown in Figure~\ref{fig:disksR}. 
In this example, the boundary components of the domain $\Omega$ surrounding the basepoint $z_0=0$ which, as expected, results in $g(r) \rightarrow 1$ more rapidly when there are more boundary components compared to when there are fewer. 

\section{Bounded domains}\label{sc:bd}
\begin{figure}[ht] %
	\centerline{\hfill
		\scalebox{0.5}{\includegraphics[trim=0 0 0 0,clip]{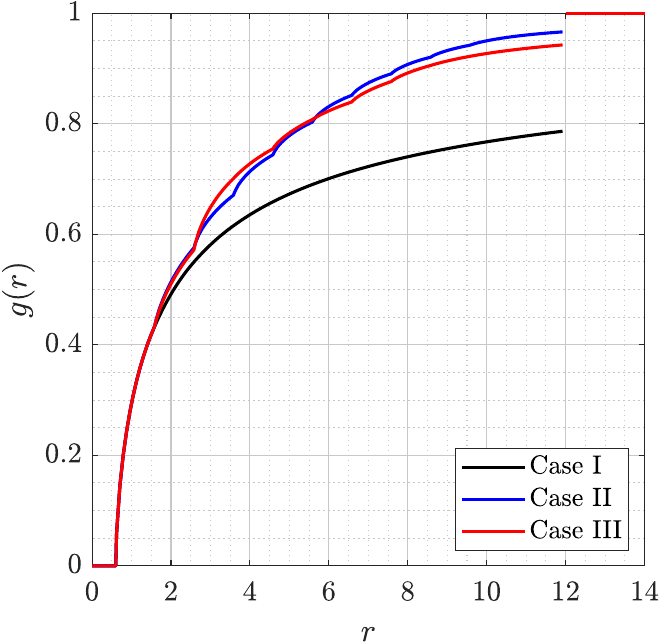}}\hfill
	}
	\caption{The values of the $g$-function $g(r)$ for the bounded domain $\Omega$ when $\ell=10$.}
	\label{fig:disksbd}
\end{figure}

To complete our study of the $g$-function in this paper, it is worth demonstrating that the presented method can be used in a straightforward manner for bounded domains $\Omega$. 
This is briefly described in this section through the following example.
Let $\Omega$ be the bounded domain lying in the interior of the circle with center $z_0=0$ and radius $R_0=12$ and exterior to the $\ell$ circles with the centers $z_k$ and radii $R_k$, $k=1,\ldots,\ell$, for the three cases considered in Section~\ref{sc:d-ex}. Then $\Omega$ is a bounded multiply connected domain of connectivity $11$.
For $0\le r<12$, the domain $\Omega_r$ here is the same as in the unbounded domain case in Figure~\ref{fig:mD10}. Hence, the values of $g(r)$ for $0\le r<12$ are the same as the values for the unbounded case as presented in Figure~\ref{fig:disks}. For all values of $r$ such that $r\ge12$, we will have $\Omega_r=\Omega$ and hence $g(r)=1$. 
The graphs of the $g$-function $g(r)$ for the three cases with $\ell=10$ are shown in Figure~\ref{fig:disksbd}. All graphs exhibit a fixed jump when $r=12$.

\section{Conclusion}\label{sc:con}

This paper has shown how to make numerical computations of the $g$-function associated with various multiply connected planar domains bounded by either rectilinear slits or circles. These $g$-functions are connected to their $h$-function counterparts, and both can be interpreted in terms of the motion of a Brownian particle released from some basepoint in the domain over which they are defined. We used a combination of conformal mapping techniques and a well-established boundary integral equation method to perform our calculations, and several graphs of $g$-functions have been plotted. 

This work is the first attempt at computing the $g$-function in multiply connected domains. The graphs of the $g$-functions presented in this paper provide some insight into the properties of these $g$-functions and how they compare to the graphs of their $h$-function counterparts shown in~\cite{green25,green24,green22}. However, much numerical work still remains to be done; in particular, to compute $g$-functions associated with domains bounded by curves other than circles and rectilinear slits, and for different basepoint locations.
Further, Stephenson's original questions in~\cite{problems} are mainly concerned about the connection between the geometry of the domain $\Omega$ and the behavior of both the $g$-function and the $h$-function. These questions are still not fully answered and hence provide a source of future research topics.

\section*{Acknowledgments} 
The authors would like to thank an anonymous reviewer for their valuable comments and suggestions. CCG acknowledges the support of a LEAPS-MPS grant from the National Science Foundation (\#2532158).

\end{document}